\documentclass[11pt]{scrartcl}
\pdfoutput=1

\usepackage{latexsym,amssymb,lastpage,natbib}
\usepackage{graphicx,amsfonts}
\usepackage{color,amsmath,fancyhdr}
\usepackage{dcolumn,textcomp,url}
\usepackage[small]{caption}
\usepackage[scaled=0.9]{beramono}
\usepackage[scaled=0.9]{helvet}
\usepackage{fourier}
\typearea{11}
\usepackage[pdftex,colorlinks]{hyperref}
\definecolor{darkblue}{cmyk}{1,0,0,0.75}
\definecolor{darkred}{cmyk}{0,1,0,0.7}
\hypersetup{anchorcolor=black, citecolor=darkblue, filecolor=darkblue, menucolor=darkblue,pagecolor=darkblue,urlcolor=darkblue,linkcolor=darkblue}
\usepackage{microtype}

\numberwithin{equation}{section}

\renewcommand{\d}{\mathop{}\!\mathrm{d}}

\begin{document}
\title{Climate tipping as a noisy bifurcation: a predictive technique}
\author{\textsc{ J. M. T. Thompson}\\[2pt]
  \begin{minipage}[t]{0.8\linewidth}
    \begin{center}
      \small Department of Applied Mathematics and Theoretical
      Physics,\\
      Centre for Mathematical Sciences, University of Cambridge\\
      Wilberforce Road, Cambridge, CB03 0WA, UK.\\
      and\\
      School of Engineering (Sixth Century Professor),
      University of Aberdeen\\[1ex]
      \
    \end{center}
  \end{minipage}\\
  \textsc{Jan Sieber}\\[2pt]    
  \begin{minipage}[t]{0.8\linewidth}
    \begin{center}
  \small
  Department of Mathematics, University of Portsmouth,\\
  Lion Gate Building, Lion Terrace,\\ Portsmouth, PO1 3HF, UK.\\[6pt]      
    \end{center}
  \end{minipage}\\
  {\small[Preprint, submitted to IMA J. Appl. Math., please cite original article]}} 

\pagestyle{fancy} 
\lhead{}
\chead{J. M. T. THOMPSON \& J. SIEBER, Climate tipping as a noisy
  bifurcation: a predictive technique}
\rhead{}
\cfoot{\thepage}
\rfoot{}
\lfoot{}
\renewcommand{\footrulewidth}{0pt}
\date{}
\maketitle


\begin{abstract} 
  \noindent {It is often known, from modelling studies, that a certain
    mode of climate tipping (of the oceanic thermohaline circulation,
    for example) is governed by an underlying fold bifurcation. For
    such a case we present a scheme of analysis that determines the
    best stochastic fit to the existing data. This provides the
    evolution rate of the effective control parameter, the variation
    of the stability coefficient, the path itself and its tipping
    point. By assessing the actual effective level of noise in the
    available time series, we are then able to make probability
    estimates of the time of tipping.  This new technique is applied,
    first, to the output of a computer simulation for the end of
    greenhouse Earth about 34 million years ago when the climate
    tipped from a tropical state into an icehouse state with ice
    caps. Second, we use the algorithms to give probabilistic tipping
    estimates for the end of the most recent glaciation of the Earth
    using actual archaeological ice-core data. }
  
\end{abstract}
\addtocontents{toc}{\protect\normalsize}
\tableofcontents

\section{Introduction}
One concern of the \emph{UN Climate Change Conference} in Copenhagen
2009 was the prediction of future climate change, subject, for
example, to a variety of carbon dioxide emission scenarios. A
particularly alarming feature of any such prediction would be a sudden
and (perhaps) irreversible abrupt change called a \emph{tipping point}
\citep{Lenton2008,Scheffer2009a}. Such events are familiar in
nonlinear dynamics, where they are called (dangerous) bifurcations, at
which one form of behaviour becomes unstable and the system jumps
rapidly to a totally different `steady state'. Many tipping points,
such as the switching on and off of ice-ages, are well documented in
paleoclimate studies over millions of years of the Earth's history.

There is currently much interest in examining climatic tipping points,
to see if it is feasible to predict them in advance using time-series
data derived from past behaviour. Assuming that tipping points may
well be governed by a bifurcation in an underlying dynamical system,
recent work looks for a slowing down of intrinsic transient responses
within the data, which is predicted to occur before most bifurcational
instabilities \citep{Held2004,Livina2007}. This is
done, for example, by determining the \emph{propagator}, which is estimated
via the correlation between successive elements of the time series, in
a window sliding along the time series. This propagator is a measure
for the linear stability. It should increase to unity at tipping.

Many trial studies have been made on climatic computer models where an
arbitrary time, $t=T$, can be chosen to represent `today'. The
challenge is then to predict a future critical time, $t^C>T$, at which
the model will exhibit a tipping instability, using only the time
history of some variable (average sea temperature, say) generated by
the model before time $T$. The accuracy of this prediction can then be
assessed by comparing it to the actual continued response of the
simulation beyond time $T$. In some cases these trials have been
reasonably successful.

Much more challenging, and potentially convincing, is to try to
predict real ancient climate tippings, using their preceding
geological data. The latter would be, for example, re-constituted
time-series provided by ice cores, sediments, etc. Using this data,
the aim would be to see to what extent the actual tipping could have
been accurately predicted in advance.

One past tipping point that has been analysed in this manner
\citep{Livina2007} is the end of the Younger Dryas event, about 11,500
years ago, when the Arctic warmed by 7\textdegree\,C in 50 years. These
authors used a time series derived from Greenland ice-core
paleo-temperature data. A second such study (one of eight made by
\citet{Dakos2008}, using data from tropical Pacific sediment cores,
gives an excellent prediction for the end of `greenhouse' Earth about
34 million years ago when the climate tipped from a tropical state
into an icehouse state.

After a review of recent research on tipping we study in
Section~\ref{sec:early:escape} the saddle-node normal form with a
drifting normal form parameter and with additive Gaussian noise to
determine the probability (or rate) of early noise-induced escape from
the potential well depending on the drift speed of the normal form
parameter and the noise amplitude.

We show how one can extract the relevant normal form quantities from a
given time series. This allows one to adjust predictions of tipping
events based on the propagator to take into account the probability of
early escape. We demonstrate how it is possible to estimate the
probability of noise-induced escape from the potential well using two
time series: one is an output from a simple stochastic model, the
other is a paleo-temperature record from ice-core data. Since the
propensity of the system to escape early from its potential well
depends only on the order of magnitude of the ratio between drift
speed and noise amplitude we expect our estimates to be reasonably
robust.

This prediction science is very young, but the above trials on
paleo-events seem very encouraging, and we describe some of them more
fully below. The prediction of \emph{future} tipping points, vital to guide
decisions about geo-engineering for example, will benefit from the
experience drawn from these trials, and will need the high quality
data currently being recorded worldwide by climate scientists today.

\section{Tipping of the Climate and its Sub-Systems}
\label{sec:tipping:climate}

\subsection{Tipping Points}
\label{sec:tipping:points}

Work at the beginning of this century which set out to define and
examine \emph{climate tipping}
\citep{Rahmstorf2001,Lockwood2001,NationalResearchCouncil2002,Alley2003,Rial2004}
focused on abrupt climate change: namely when the Earth system is
forced to cross some threshold, triggering a sudden transition to a
new state at a rate determined by the climate system itself and
(usually) faster than the cause, with some degree of
irreversibility. Recently, the \emph{Intergovernmental Panel on
  Climate Change} \citep{IPCC2007} made some brief remarks about
abrupt and rapid climate change, while \citet{Lenton2008} have sought
to define these points more rigorously. The physical mechanisms
underlying these tipping points are typically internal positive
feedback effects of the climate system (thus, a certain propensity for
saddle-node bifurcations).

\subsection{Tipping Elements}
\label{sec:tipping:elements}

In principle a climate tipping point might involve simultaneously many
features of the Earth system, but it seems that many tipping points
might be strongly associated with just one fairly well defined
sub-system. These \emph{tipping elements} are well-defined sub-systems
of the climate which work (or can be assumed to work) fairly
independently, and are prone to sudden change. In modelling them,
their interactions with the rest of the climate system are typically
expressed as a control parameter (or forcing) that varies slowly over
time.

Recently \citet{Lenton2008} have listed nine tipping elements that
they consider to be primary candidates for future tipping due to human
activities, and as such have relevance to political decision making at
\citet{Copenhagen2009} and beyond. These elements, their possible outcomes,
and Lenton's assessment of whether their tipping might be associated
with an underlying bifurcation are:
\begin{itemize}
\item[(1)] loss of Arctic summer sea-ice
(possible bifurcation);
\item[(2)] collapse of the Greenland ice sheet
(bifurcation);
\item[(3)] loss of the West Antarctic ice sheet (possible
bifurcation);
\item[(4)] shut-down of the Atlantic thermohaline circulation
(fold bifurcation);
\item[(5)] increased amplitude or frequency of the El
Ni{\~n}o Southern Oscillation (some possibility of bifurcation);
\item[(6)] switch-off of the Indian summer monsoon (possible
  bifurcation);
\item[(7)] changes to the Sahara/Sahel and West African
  monsoon, perhaps greening the desert (possible bifurcation);
\item[(8)] loss of the Amazon rainforest (possible bifurcation);
\item[(9)] large-scale dieback of the Northern Boreal forest (probably
  not a bifurcation).
\end{itemize}
The analysis and prediction of tipping points of climate subsystems is
currently being pursued in several streams of research, and we should
note in particular the excellent book by Marten Scheffer about tipping
points in `Nature and Society', which includes ecology and some
climate studies \citep{Scheffer2009a}.
\section{Bifurcations and their Precursors}
\label{sec:bif:prec}

\subsection{Generic Bifurcations of Dissipative Systems}
\label{sec:genbif}
The great revolution of nonlinear dynamics over recent decades has
provided a wealth of information about the bifurcations that can
destabilise a slowly evolving system like the Earth's climate. These
bifurcations are defined as points during the slow variation of a
`control' parameter at which a qualitative topological change of
behaviour is observed in the multi-dimensional phase space of the
system.

The Earth's climate is what dynamicists would call a \emph{dissipative
  system}, and for this the bifurcations that can be typically
encountered under the variation of a single control parameter are
classified into three types, \emph{safe}, \emph{explosive} and
\emph{dangerous} \citep{Thompson1994,Thompson2002}.

The safe bifurcations, such as the supercritical Hopf bifurcation,
exhibit a continuous supercritical growth of a new attractor path with
no fast jump or enlargement of the attracting set. They are
determinate with a single outcome even in the presence of small noise,
and generate no hysteresis with the path retraced on reversal of the
control sweep. The explosive bifurcations are less common phenomena
lying intermediate between the safe and dangerous types: we simply
note here that, like the safe bifurcations, they do not generate any
hysteresis. The dangerous bifurcations are typified by the simple fold
(saddle-node bifurcation) at which a stable path increasing with a
control parameter becomes unstable as it curves back towards lower
values of the control, and by the subcritical bifurcations. They
exhibit the sudden disappearance of the current attractor, with a
consequential sudden jump to a new attractor (of any type). They can
be indeterminate in outcome, depending on the topology of the phase
space, and they always generate hysteresis with the original path not
reinstated on control reversal.

Any of these three bifurcation types could in principle underlie a
climate tipping point. But it is the dangerous bifurcations that will
be of major concern, giving as they do a sudden jump to a different
steady state with hysteresis, so that the original steady state will
not be re-instated even if the controlling cause is itself
reversed. So any future climatic tipping to a warmer steady state may
be irreversible: a subsequent reduction in CO$_2$ concentration will not
(immediately, or perhaps ever) restore the system to its pre-tipping
condition.

\subsection{Time-Series Analysis of Incipient Bifurcations}
\label{sec:timeseries}
Most of the bifurcations in dissipative systems, including the static
and cyclic folds which are the most likely to be encountered in
climate studies, have the following useful precursor (see
\citet{Thompson2002} for more details). The stability and attracting
strength of the current steady state is becoming steadily weaker and
weaker in one mode as the bifurcation point is approached. This
implies that under inevitable noisy disturbances, transient motions
returning to the attractor will become slower and slower: in the
limit, the rate of decay of the transients decreases linearly to zero
along the path. 

\citet{Held2004} and \citet{Livina2007} have recently presented
algorithms that are able to detect incipient saddle-node bifurcations
from time series of dynamical systems. Both methods estimate the
linear decay rate (LDR) toward a quasi-stationary equilibrium that is
assumed to exist and to drift toward a saddle-node bifurcation.
Typical test data for the algorithms comes from either archaeological
records or from output of climate models.  Both algorithms
(\emph{degenerate finger-printing} by \citet{Held2004} and
\emph{detrended fluctuation analysis} by \citet{Livina2007}) have to
make assumptions about the process underlying the recorded time series
that are generally believed to be sensible for the tipping elements
listed by \citet{Lenton2008}.  First, one has to assume that the
process is a dynamical system close to a stable equilibrium that
drifts only slowly but is perturbed by (random) disturbances. The
second assumption is that the system is effectively one-dimensional,
that is, the equilibrium of the undisturbed system is strongly stable
in all directions except a single critical one. Quantitatively this
means that one assumes the presence of three well-separated time
scales, expressed as rates:
\begin{displaymath}
  \kappa_\mathrm{drift}\ll\kappa_\mathrm{crit}\ll\kappa_\mathrm{stab}\mbox{.}
\end{displaymath}
Here $\kappa_\mathrm{drift}$ is the average drift rate of those
quantities that the algorithm treats as a parameter, for example,
freshwater forcing in studies of the thermohaline circulation (THC,
the global heat- and salinity-driven conveyor belt of oceanic
water). The rate $\kappa_\mathrm{crit}$ is the rate with which a small
disturbance in state space relaxes back to equilibrium. The rate
$\kappa_\mathrm{stab}$ is the decay rate of all other
\emph{non-critical} modes. We note that the drift rate
$\kappa_\mathrm{drift}$ becomes \emph{larger} than
$\kappa_\mathrm{crit}$ once the drifting parameter is very close to
its bifurcation value. Third, one assumes that the disturbances are
small in the sense that the relaxation to equilibrium is governed
mostly by the linear decay rates (this implies, for example, that the
potential well in which the dynamical system can be imagined to be
sitting is approximately symmetric in the critical direction).

The basic procedure proposed by \citet{Held2004} consists of
three steps, given a time series $(t_k,z_k)$ of measurements $z_k$ at
possibly unevenly spaced time points $t_k$.
\begin{enumerate}
\item \textbf{Interpolation} Choose a stepsize $\Delta t$ satisfying
  \begin{displaymath}
    \kappa_\mathrm{crit}\ll\frac{1}{\Delta t}\ll\kappa_\mathrm{stab}
  \end{displaymath}
  and interpolate such that the spacing in time is uniform.  Now, one
  has a new time series
  \begin{displaymath}
    z_{k,\mathrm{new}}=z(k \Delta t)\mbox{,}
  \end{displaymath}
  evenly spaced in time.
\item \textbf{Detrending} Remove the slow drift of the equilibrium by
  subtracting a slowly moving average. For example, choose a Gaussian
  kernel
  \begin{displaymath}
    G_k(t)=\frac{1}{\sqrt{2\pi}d}\exp\left(-\frac{1}{2}\frac{(t-k\Delta t)^2}{d^2}\right)
  \end{displaymath}
of bandwidth $d$ satisfying
  \begin{displaymath}
    \kappa_\mathrm{drift}\ll\frac{1}{d}\ll\kappa_\mathrm{crit}\mbox{,}
  \end{displaymath}
  and subtract the average
  \begin{displaymath}
    Z(k\Delta t)=\frac{\sum_{i=1}^N
      G_k(i\Delta t)z_i}{\sum_{i=1}^NG_k(i\Delta t)}
  \end{displaymath}
  of $z_k$ over the kernel. The result of this is a time series
  \begin{displaymath}
    y_k=z_k-Z(k\Delta t)
  \end{displaymath}
  which fluctuates around zero and can be considered as stationary on
  time scales shorter than $1/\kappa_\mathrm{drift}$.
\item \textbf{Fit LDR in sliding window} One assumes that the
  remaining time series, $y_k$, can be modelled by a stable scalar linear
  mapping disturbed by noise, a so-called AR(1)-model
  \begin{equation}\label{eq:ar1}
    y_{k+1}=c_k y_k   + \theta\eta_k
  \end{equation}
  where $\theta\eta_k$ is the instance of a random disturbance of
  amplitude $\theta$ at time $k\Delta t$ and $c_k$ is the
  \emph{propagator}, related to $\kappa_\mathrm{crit}$ at time
  $k\Delta t$ via
  \begin{displaymath}
    c_k=\exp(-\kappa_{\mathrm{crit},k}\Delta t)\mbox{.}    
  \end{displaymath}
  If one assumes that the disturbances $\eta_k$ have a normal
  distribution and are independent from each other, and that $c_k$ is
  nearly constant on time scales shorter than
  $1/\kappa_\mathrm{drift}$, one can choose a sliding window size
  $w=2m+1$ satisfying
  \begin{displaymath}
    \kappa_\mathrm{drift}\ll \frac{1}{w\Delta t}    
  \end{displaymath}
 and determine the
  propagator $c_k$ by an ordinary least-squares fit of
  \begin{displaymath}
    y_{j+1}=c_k y_j    
  \end{displaymath}
  over the set of indices $j=k-m\ldots k+m$. An estimate for the noise
  amplitude $\theta$ can be obtained from the standard deviation of
  the residual of the linear least-squares fit:
  \begin{displaymath}
    \theta_k=\operatorname{stdev}%
    \left([y_j-c_ky_j]_{j=k-m}^{\phantom{j=}k+m}\right)
  \end{displaymath}
\end{enumerate}

\begin{figure}[ht!]
  \centering
  \includegraphics[scale=0.6]{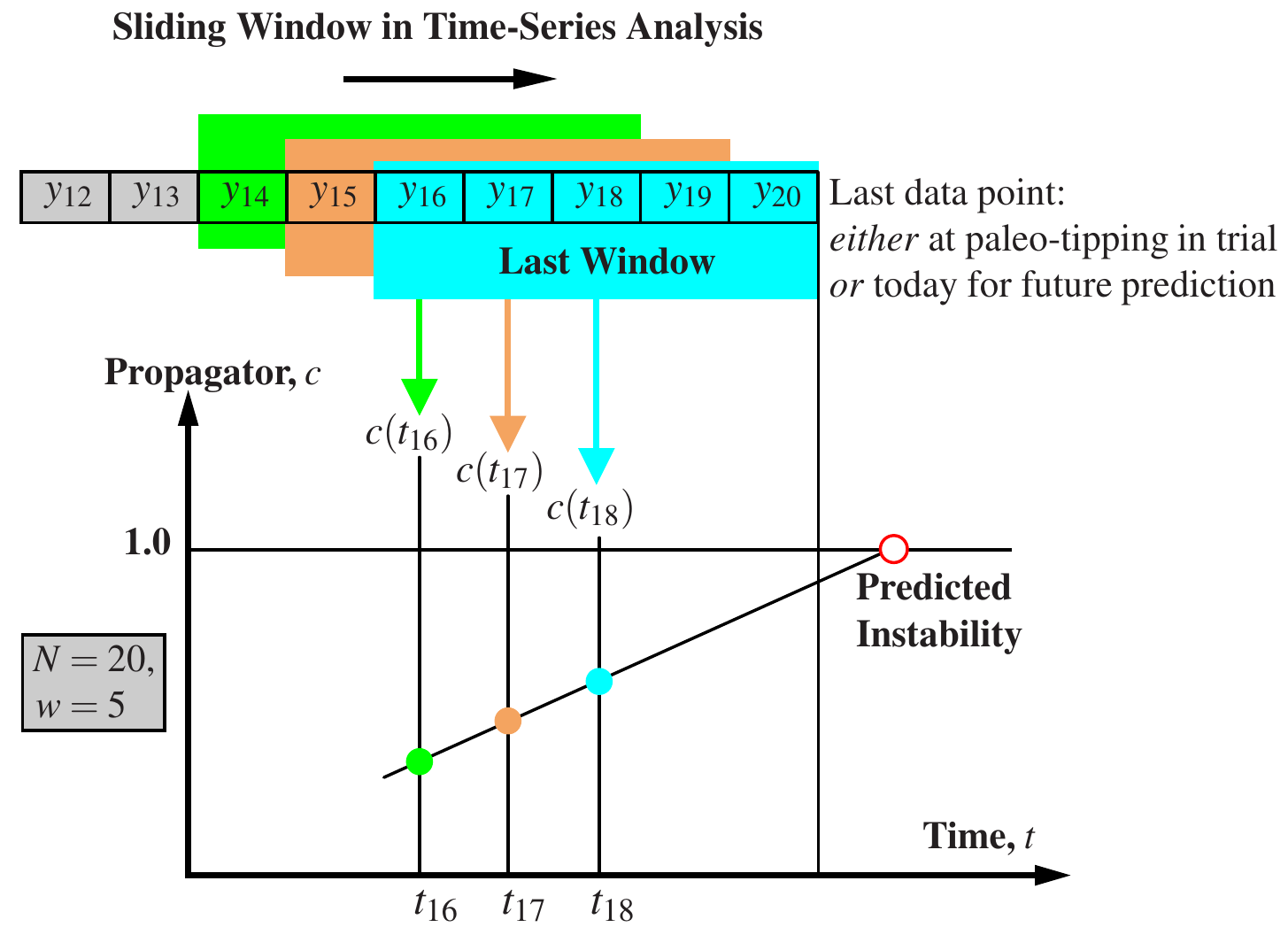}
  \caption{Illustration of sliding windows used in degenerate
    finger-printing proposed by \citet{Held2004}. The extracted
    autocorrelation coefficient $c_k$ is an estimate for the
    propagator at the mid-time of the sliding windows.  Extrapolation
    is required for prediction.}
  \label{fig:sliding}
\end{figure}

This process will stop when the front end of the last window hits the
last data point. Notice that when using paleo-data the end of the
analysed time series should be chosen before the tipping point that is
the object of the investigation. This choice is essential to prevent
data, spurious to our predictions of the pre-tipping behaviour,
entering the analysis: it also makes a complete analogy with any
attempt to predict future tipping points from data terminating
today. Figure~\ref{fig:sliding} illustrates this requirement on the
sliding windows (with a time series of length $N=20$ for
illustration). It also shows a simple linear extrapolation that one
might make in order to predict where tipping occurs.
In this manner \citet{Held2004} obtain the so-called propagator graph
of the estimated $c_k$ versus its central time, as illustrated in
Figure~\ref{fig:sliding}. On this graph, $c$ is expected to head
towards $+1$ at any incipient bifurcation. In other words, the slowing
down of the relaxation from disturbances along the times series can
serve as an early-warning signal for an imminent bifurcation
\citep{Dakos2008}. We should finally note that having used a
first-order mapping in Equation~\eqref{eq:ar1} and employed
autocorrelation techniques, the propagator $c$ is often called the
first-order autoregressive coefficient and written as ARC(1). The
prediction based on detrended fluctuation analysis, as proposed by
\citet{Livina2007}, also reconstructs the propagator $c$ but does so
via the scaling exponent of the variance of the (detrended) time
series. For a more complete description of the time-series techniques
employed see the recent review by \citet{Thompson2010}. Also, see
\citet{Corsi2007} for a discussion of similar problems in power systems
engineering (the prediction of fold bifurcations leading to voltage
collapse in energy networks).

\section{Review of Recent Work}
\label{sec:review}

\subsection{First Prediction of an Ancient Tipping}
\label{sec:pred}
The first prediction of an ancient climate tipping event using
preceding geological data is due to \citet{Livina2007} who tested
their detrended fluctuation analysis on the rapid warming of the earth
that occurred about 11,500 years ago at the end of the so-called
Younger Dryas event (analysing Greenland ice-core paleo-temperature
data, which is available from 50,000 years ago to the present).

This Younger Dryas event \citep{Houghton2004} was a curious cooling
just as the Earth was warming up after the last ice age, as is clearly
visible, for example, in records of the oxygen isotope $\delta^{18}O$
in Greenland ice. It ended in a dramatic tipping point, 11,500 years
ago, when the Arctic warmed by 7\textdegree\,C in 50 years. Its
behaviour is thought to be linked to changes in the thermohaline
circulation (THC). This `conveyor belt' is driven by the sinking of
cold salty water in the North and can be stopped if too much
fresh-melt makes the water less salty, and so less dense. At the end
of the ice age when the ice-sheet over North America began to melt,
the water first drained down the Mississippi basin into the Gulf of
Mexico. Then, suddenly, it cut a new channel near the St Lawrence
river to the North Atlantic. This sudden influx of fresh water cut off
part of the ocean `conveyor belt', the warm Atlantic water stopped
flowing North, and the Younger Dryas cooling was started. It was the
\emph{re-start} of the circulation that could have ended the Younger
Dryas at its rapid tipping point, propelling the Earth into the warmer
Pre-Boreal era.

\begin{figure}[ht!]
  \centering
  \includegraphics[scale=0.6]{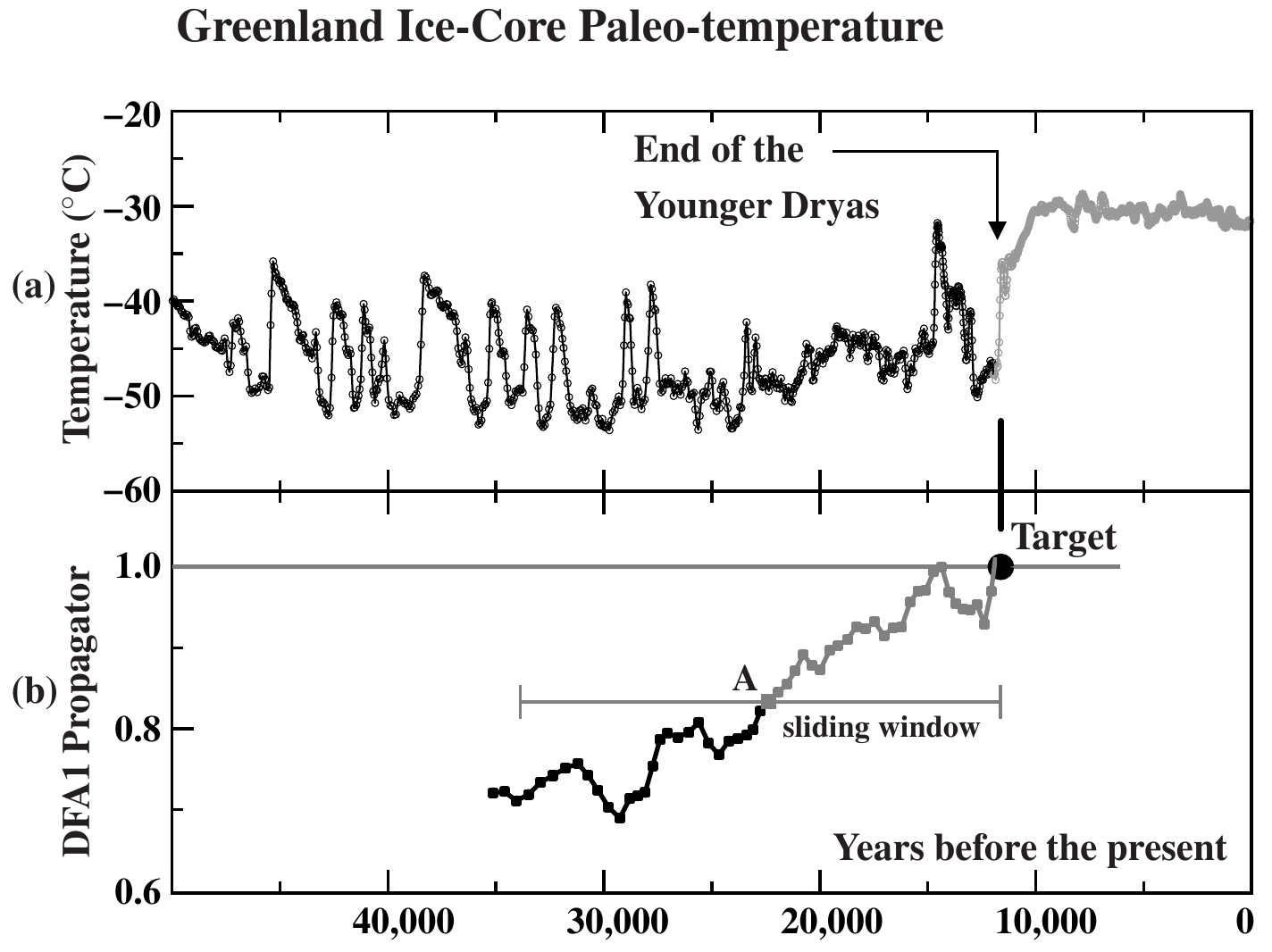}
  \caption{Results of \citet{Livina2007}: (a) Greenland ice-core
    (GISP2) paleo-temperature with an unevenly spaced record, visible
    in the varying density of symbols on the curve. The total number
    of data points is $N=1586$. In (b) the DFA1-propagator is
    calculated in sliding windows of length $w=500$ points and mapped
    into the middle points of the windows. A typical sliding window
    ending near the tipping is shown.}
  \label{fig:livina}
\end{figure}

The results of \citet{Livina2007} are shown in
Figure~\ref{fig:livina}, where their propagator (based on detrended
fluctuation analysis, DFA) is seen heading towards its critical value
of $+1$ at about the correct time. Notice, though, that from a
prediction point of view the propagator graph should end at point A
when the estimation-window reaches the tipping point. In this example
extracting the propagator is particularly challenging because the data
set was comparatively small (1586 points) and unevenly spaced.
\subsection{Systematic Study of Eight Ancient Tippings}
\label{sec:dakos}
In a more recent paper, \citet{Dakos2008} systematically estimated a
propagator stability coefficient from reconstructed time series of
real paleo-data preceding eight ancient tipping events. These are:
\begin{itemize}
\item[(a)] the end of the greenhouse Earth about 34 million years ago
  when the climate tipped from a tropical state (which had existed for
  hundreds of millions of years) into an icehouse state with ice caps,
  using data from tropical Pacific sediment cores;
\item[(b)] the end of the last glaciation, and the ends of three
  earlier glaciations, drawing data from the Antarctica Vostok ice
  core;
\item[(c)] the B{\o}lling-Aller{\"o}d transition which was dated about
  14,000 years ago, using data from the Greenland GISP2 ice core;
\item[(d)]
  the end of the Younger Dryas event about 11,500 years ago, as
  discussed in the previous section, but drawing not on the Greenland
  ice core, but rather on data from the sediment of the Cariaco basin
  in Venezuela;
\item[(e)] the desertification of North Africa when there was a sudden
  shift from a savanna-like state with scattered lakes to a desert
  about 5,000 years ago, using the sediment core from ODP Hole 658C,
  off the west coast of Africa.
\end{itemize}
In all of the cases studied by \citet{Dakos2008}, the propagator $c_k$
as extracted by degenerate finger-printing was shown to exhibit a
statistically significant increase (corresponding to a slowing down of
the relaxation) prior to the tipping transition. \citet{Dakos2008}
also demonstrated that their principal result, the statistically
significant increase of $c_k$, is robust with respect to variations in
smoothing kernel bandwidth $d$, sliding window length $w$ and the
interpolation procedure.

\section{Noise-induced systematic bias of extrapolated prediction}
\label{sec:early:escape}
The intention behind the development of the time series analysis
algorithms goes beyond statistical evidence of an increasing LDR: the
goal of both algorithms is to \emph{predict the time} (or probability)
of the tipping event from the observational data before the event
takes place. This is more challenging and suffers from additional
uncertainties. Apart from the dependence of the value of $c_k$ on
algorithm parameters (for example the sliding window length $w$), for
a prediction of the time of tipping we have to extrapolate. This
implies that we have to assume that the underlying control parameter
drifts with nearly constant speed during the recorded time
series. This is often not the case in the study of subsystems of the
climate when the control parameter is determined by the dynamics of
another, coupled, subsystem. Even if the control parameter drifts with
constant speed, for prediction we have to assume in addition
\emph{transversality}, that is, the control parameter has to vary the
unfolding parameter of the normal form of the saddle-node bifurcation
nearly linearly.

\begin{figure}[!t]
  \centering
  \includegraphics[scale=0.55]{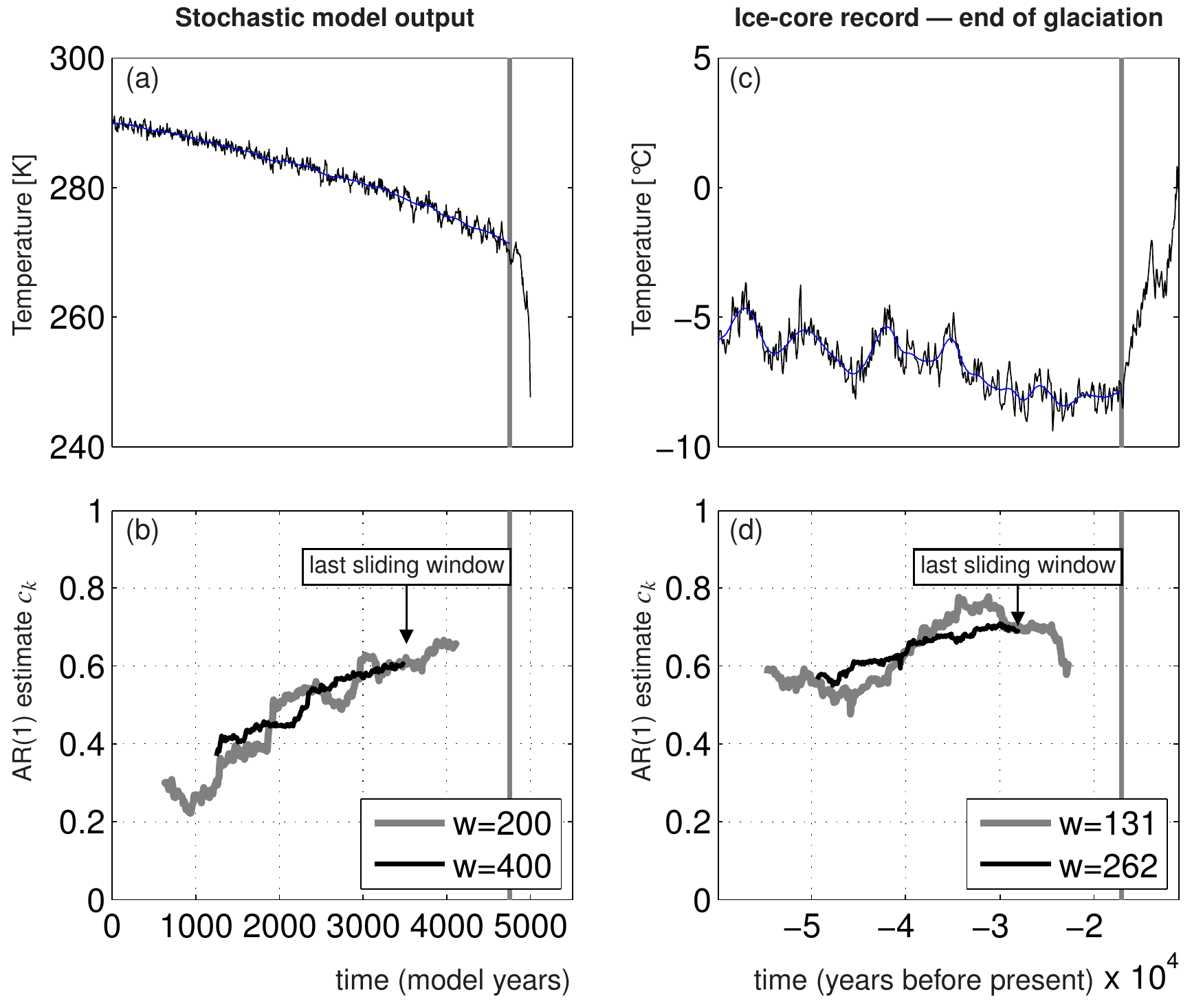}
  \caption{LDR estimates for the output of a model of transition to
    icehouse earth ((a) and (b)), and for an archaeological
    temperature record of the end of the recent glaciation ((c) and
    (d)). See \citep{Petit1999} and main text for description. Data
    source:
    \protect\url{ftp://ftp.ncdc.noaa.gov/pub/data/paleo/icecore/antarctica/vostok/deutnat.txt}.}
  \label{fig:ldr:est}
\end{figure}
Figure~\ref{fig:ldr:est} shows two time series ((a) and (c)) and the
corresponding time series of extracted estimates for the propagator
$c_k$ ((b) and (d)). Time series (a) is the output of a model
simulation for a transition to an icehouse Earth, and is taken from
\citep{Dakos2008}. The model as presented by \citet{Dakos2008} is a
scalar stochastic ODE where a control parameter is varied linearly in
time and the system is known to encounter a saddle-node bifurcation
(originally the model was developed by \citet{Fraedrich1978}; see
supplement of \citep{Dakos2008}). Time series (b) shows the propagator
$c_k$ extracted from time series (a) using degenerate
finger-printing. Time series (c) is a snapshot of temperatures before
the end of the last glaciation, 20,000 years ago. The data is taken
from \citep{Petit1999}, the window of the snapshot is identical to
Figure 1(I) in \citep{Dakos2008}. The estimated propagator $c_k$ for
time series (c) is shown in diagram (d). The graphs in diagrams (b)
and (d) are qualitatively the same as in \citep{Dakos2008} but differ
slightly quantitatively. This is likely due to minor variations
between the algorithm parameters that we used and the algorithm
parameters used by \citet{Dakos2008}. To provide a visual clue about
the level of uncertainty in the time series of $c_k$ the diagrams (b)
and (d) show estimates extracted using two different window size
parameters $w$. Despite the uncertainty two features of the time
series of propagators $c_k$ are discernible: first, as studied in
detail in \citep{Dakos2008}, $c_k$ is increasing. Second, linear
extrapolation does not match if we expect the tipping to occur at the
extrapolated time for $c=1$ (which would correspond to the critical
value of the propagator). In both cases the tipping occurs earlier,
and at least for the model output (which is known to drift through a
saddle-node bifurcation) the bias of the linear extrapolation is
systematic. There are two competing effects which might determine the
systematic bias. On the one hand the equilibrium starts to drift
rapidly (square-root like) when the drifting parameter approaches the
saddle-node bifurcation, on the other hand random disturbances may
kick the dynamical system out of the shrinking basin of attraction
prematurely. The numerical study in the following sections quantifies
these two effects.

\subsection{Rate of noise-induced escape from basin near saddle-node}
\label{sec:snform}
We take the normal form of a saddle-node bifurcation, perturb it by
adding Gaussian white noise and let the control parameter $a$ drift
with speed $\epsilon$ (the caricature climate model of
\citet{Fraedrich1978} is exactly of this type):
\begin{align}
  \label{eq:gensnform}
  \d x&=\left[a-x^2\right]\d t+\sigma\d W_t\\
  \label{eq:adrift}
  \d a&=-\epsilon \d t
\end{align}
The perturbation $W_t$ (sometimes written $W(t)$ to avoid double
subscripts) is a Wiener process, which is defined by two properties:
\begin{enumerate}
\item $W_0=0$, and
\item for every sequence of time points $0\leq t_1 \leq t_2 \leq
  \ldots \leq t_k$ all increments $W(t_{i+1})-W(t_i)$ are
  independent random variables with Gaussian distribution of zero mean
  and variance $t_{i+1}-t_i$ (and, thus, standard deviation
  $\sqrt{\smash[b]{t_{i+1}-t_i}}$).
\end{enumerate}
The coefficient $\sigma$ controls the amplitude of the noise variance
added to the slowly drifting equilibrium. If we \emph{freeze} the
drifting (set $\epsilon=0$) and set the noise amplitude to zero
($\sigma=0$) then the dynamics of \eqref{eq:gensnform} corresponds to
the dynamics of an overdamped particle in a potential well of the
shape
\begin{displaymath}
U(x)=-ax+x^3/3\mbox{.}  
\end{displaymath} 
This dynamics of Equation~\eqref{eq:gensnform} with a fixed $a$ and no
noise has a stable equilibrium $X_s$ at $\sqrt{\smash[b]{a}}$ and an
unstable equilibrium $X_u$ at $-\sqrt{\smash[b]{a}}$.
Correspondingly, the potential well $U$ has a (local) minimum at
$X_s=\sqrt{\smash[b]{a}}$ and a hill-top (local maximum) at
$X_u=-\sqrt{\smash[b]{a}}$. For $x$ going to $-\infty$ the potential
$U$ falls off to $-\infty$, and for $x$ going to $+\infty$ it
increases to $+\infty$.  The differential equation \eqref{eq:adrift}
for $a$ governs the drifting of the control and has the solution
$a(t)=a(0)-\epsilon t$.

We say that a trajectory $x(t)$ \emph{escapes} if it reaches $-\infty$
(typically in finite time). In numerical tests one can detect if the
trajectory crosses a fixed negative threshold $x_\mathrm{th}$ from
which it is is unlikely to return to the well (for example, set
$x_\mathrm{th}$ at some fixed distance below $X_u$). We are interested
in finding the cumulative escape probability $P_\mathrm{esc}(a)$ for a
trajectory to escape before the drifting control parameter has reached
the value $a$. If the stochastic process
\eqref{eq:gensnform}--\eqref{eq:adrift} starts with a sufficiently
large $a(0)$ and if $\epsilon$ is sufficiently small then this
probability $P_\mathrm{esc}$ depends only weakly on the initial
distribution of $x$ as long as this initial distribution is
concentrated inside the potential well.

The perturbation $W_t$ is a random disturbance (noise) as if $x$ is
the position of a particle that has a temperature. If the modulus of
the noise is orders of magnitude smaller than the height of the
barrier (the potential difference between $X_s$ and $X_u$, which is
$4/3\cdot a^{3/2}$\,) then the escape rate formulas developed for
chemical reaction rates (Kramers' escape rate, see \citep{Hanggi1990})
can be applied. As we are interested in the transition of $a$ through
zero the standard reaction rate theory is not applicable. Analytical
formulas for jumping times in periodic potentials from one period to
the next and arbitrary noise are given in \citep{Malakhov1996a}. See
also \citep{Fischer2006} for an analysis of stochastic resonance
between two slowly varying potential wells.

\paragraph*{Re-scaling of parameters}
In the noise-perturbed drifting normal form
\eqref{eq:gensnform}--\eqref{eq:adrift} we can re-scale time $t$, $x$,
$a$ and $\epsilon$ such that the noise amplitude $\sigma$ is equal to
$1$. Exploiting that $\d W_{\rho s}=\sqrt{\smash[b]{\rho}}\d W_s$, and
setting
\begin{equation}
  \label{eq:rescaling}
  t_\mathrm{old}\sigma^{2/3}=t_\mathrm{new}\mbox{,\quad}
  \sigma^{-2/3}x_\mathrm{old}=x_\mathrm{new}\mbox{,\quad}
  \sigma^{-4/3}a_\mathrm{old}=a_\mathrm{new}\mbox{, and\quad}
  \sigma^{-2}\epsilon_\mathrm{old}=\epsilon_\mathrm{new}\mbox{,} 
\end{equation}
we obtain a stochastic process for the re-scaled quantities that is of
the same form as the original noise-perturbed normal form
\eqref{eq:gensnform}--\eqref{eq:adrift} (except that the parameter
$\sigma$ has been absorbed as a unit):
\begin{align}
  \label{eq:snformscal}
  \d x&=\left[a-x^2\right]\d t+\d W_t\\
  \label{eq:adriftscal}
  \d a&=-\epsilon \d t\mbox{.}
\end{align}
In the rescaled coordinates the node (the local minimum of the well)
and the saddle (the barrier) are at
\begin{displaymath}
  X_{u,\mathrm{new}}=\sigma^{-2/3}X_{u,\mathrm{old}}\mbox{\quad and\quad}
  X_{s,\mathrm{new}}=\sigma^{-2/3}X_{s,\mathrm{old}}\mbox{,}  
\end{displaymath}
respectively.

\paragraph*{Escape rates for the saddle-node normal form with noise}
For sufficiently small drift speeds $\epsilon$ of the control
parameter $a$ we can approximate the escape probability
$P_\mathrm{esc}(a)$ in the drifting system
\eqref{eq:snformscal}--\eqref{eq:adriftscal} using quantities of the
frozen system (drift speed $\epsilon=0$). A useful quantity is the
\emph{escape rate} $k_0$. For fixed $a>0$ in the noise-perturbed
normal form \eqref{eq:snformscal} the escape rate $k_0(a)$ can be
defined as follows:
\begin{enumerate}
\item set a large ensemble (size $N$) of initial values $x$ inside the
  potential well (for example, $x=X_s$);
\item evolve the noise-perturbed normal form \eqref{eq:snformscal} and
  measure the fraction $r(t)$ of instances that have reached $-\infty$
  during the previous unit time interval:
  \begin{displaymath}
    r(t)=\lim_{h\to0}\,\lim_{N\to\infty}\frac{1}{h}\, 
    \frac{\mbox{number of instances reaching $-\infty$ 
        during $[t,t+h]$}}{\mbox{number of instances still 
        finite at time $t$}}
  \end{displaymath}
\item after a transient this fraction $r(t)$ converges to a constant
  $k_0$.
\end{enumerate}
This convergence is achieved only in the limit $N\to\infty$ due to
depletion for finite ensembles. In numerical calculations one can
replace $-\infty$ with a finite threshold $x_\mathrm{th}\ll X_u$. One
can avoid depletion of the finite ensemble by re-initialising the
escaped instance to an initial value that is randomly selected from
the remaining (non-escaped) ensemble (``putting the particle back into
the well'').

\begin{figure}[!t]
  \centering
  \includegraphics[width=0.8\textwidth]{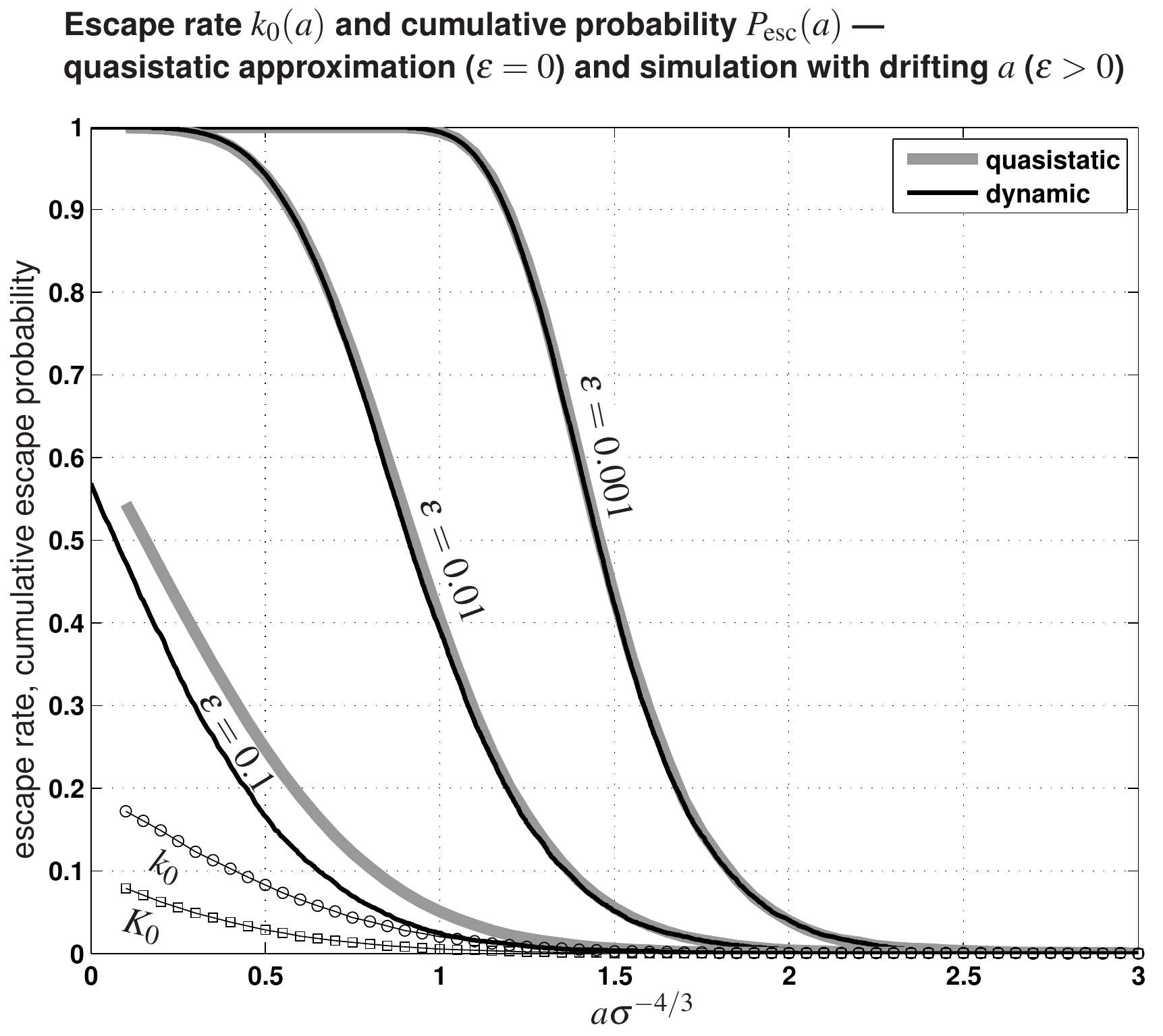}
  \caption{Escape rate $k_0$ from potential well for fixed $a$ (in the
    original parameters of
    \eqref{eq:gensnform}--\eqref{eq:adrift}). The function $K_0(a)$ is
    the integral of $k_0$, as defined in \eqref{eq:probesc},
    integrated from $a$ to $3$. The curves $\epsilon=0.1$,
    $\epsilon=0.01$, $\epsilon=0.001$ are the cumulative probabilities
    for escape from the potential well for dynamically decreasing
    depth of the potential well before $a$ is reached. The quantity
    $\epsilon$ indicates the rate of change of the potential in the
    unscaled parameters. For comparison the quasi-static estimate for
    the cumulative probability of escape is indicated using grey
    curves (see also Figure~\ref{fig:percentiles}).}
  \label{fig:snfprob}
\end{figure}
Figure~\ref{fig:snfprob} shows this escape rate $k_0$ as a function of
$a$ as a grey curve with circles. Note that the $x$-axis corresponds
to the re-scaled $a$ after transformation \eqref{eq:rescaling}. An
escape rate $k_0\approx0.1$ at $a\sigma^{-4/3}=0.4$ means that
approximately $10\%$ of the realizations cross the threshold
$x_\mathrm{th}=-5$ to escape per unit time interval during the
solution of the saddle-node normal form \eqref{eq:gensnform} with
$a\sigma^{-4/3}=0.4$.

\paragraph*{Quasi-static approximation of the cumulative escape probability}
\begin{figure}[ht]
  \centering
  \includegraphics[scale=0.5]{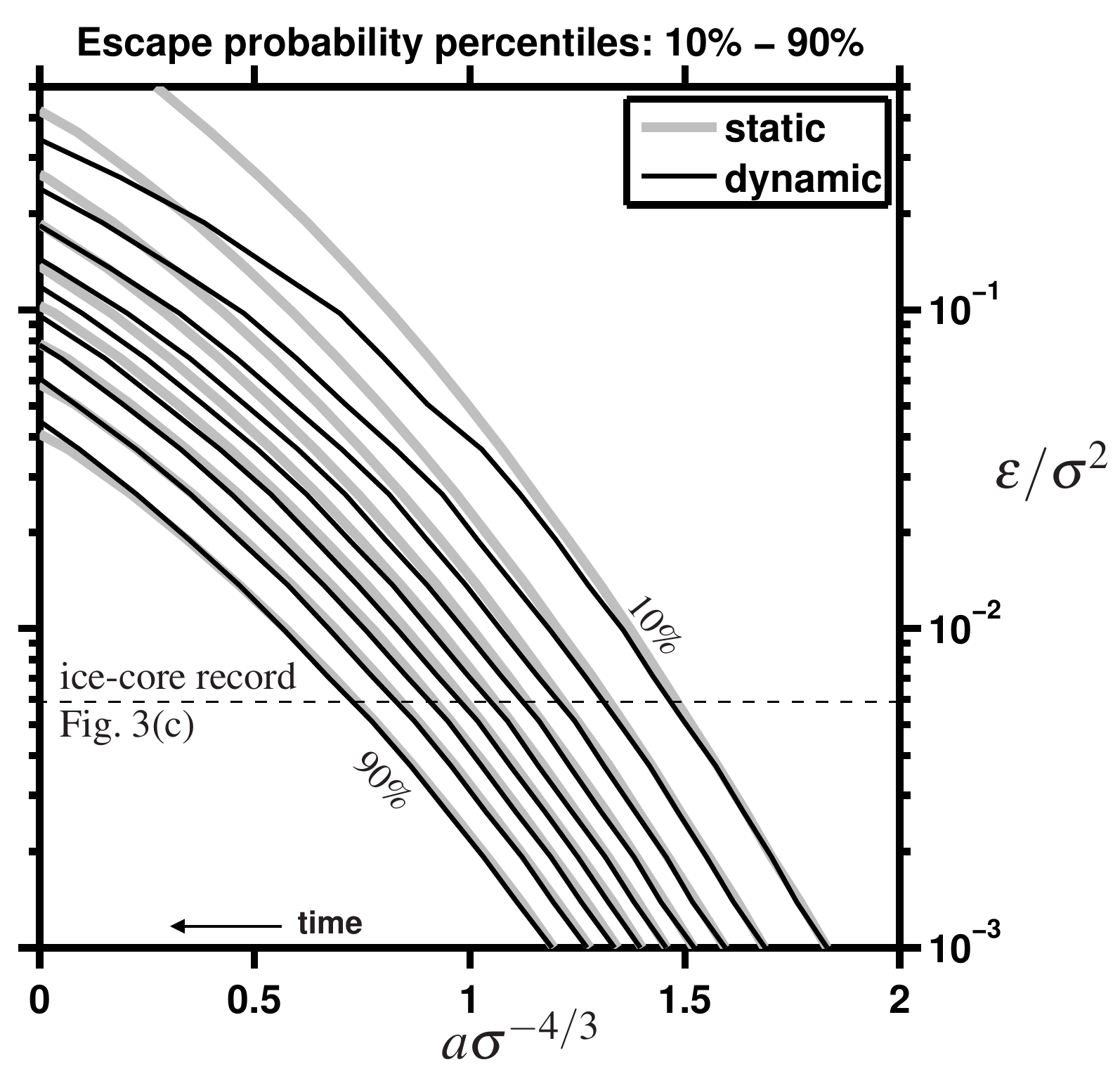}
  \caption{Percentiles of the probability of escape for dynamically
    varying $a$ (in black) in the original parameters of
    \eqref{eq:gensnform}--\eqref{eq:adrift}. For fixed $\epsilon$ the
    $x$\%-curve shows at which $a$ $x$\% of the realisations have
    escaped. For comparison the grey curves show the estimate obtained
    using the quasi-static approximation. The depth of the potential
    well (and, thus, $a$) decreases in time.}
  \label{fig:percentiles}
\end{figure}
Using the escape rate $k_0(a)$ for the frozen problem ($\epsilon=0$ in
\eqref{eq:snformscal}--\eqref{eq:adriftscal}) we can approximate the
cumulative probability $P_\mathrm{esc}$ of escape for a trajectory of
the saddle-node normal form with drifting control parameter. If we
assume that escape is irreversible and denote by $p(t)$ the
probability that the trajectory has not escaped until time $t$ then we
have the relation
\begin{equation}\label{eq:condprob}
  p(t+h)=\left[1-(k_0(a(t))+O(\epsilon))h\right]p(t)
\end{equation}
for small time steps $h>0$. The factor $1-k_0(a(t))h$ is the
probability that the trajectory escapes during the time interval
$[t,t+h]$ if we approximate the slowly changing variable $a$ by its
left-end value $a(t)$ in this interval. Relation~\eqref{eq:condprob}
expresses that a trajectory will not escape until time $t+h$ if it has
not escaped until time $t$ and it does not escape during the interval
$[t,t+h]$. Letting $h$ go to zero and ignoring the slow drift of $a$
during the time interval $[t,t+h]$ of order $O(\epsilon)$ we obtain
the differential equation
\begin{displaymath}
  \frac{\d}{\d t}p(t)=-k_0(a(t))p(t)\mbox{,}
\end{displaymath}
which has the solution
\begin{equation}\label{eq:probstayt}
  p(t)=\exp\left(\int_0^t-k_0(a(s))\d s\right)
\end{equation}
if we start with an initial distribution concentrated in the potential
well ($p(0)=1$). We want to find the approximate cumulative
probability $P_\mathrm{esc}(a)$ of escape before the control 
has drifted to a certain value $a$, so we substitute
$a=a(0)-\epsilon t$ into expression \eqref{eq:probstayt} for $p(t)$:
\begin{equation}
  \label{eq:probesc}
  P_\mathrm{esc}(a)=
  1-\exp\left(-K_0(a)/\epsilon\right)
  \mbox{\quad where \quad}
  K_0(a)=\int_{a}^{a(0)}k_0(a')\d a'\mbox{.}
\end{equation}
This approximation \eqref{eq:probesc} for the escape probability
assumes that $\epsilon$ is small, that the escape is irreversible
(which is accurate if the threshold $x_\mathrm{th}$ is sufficiently
negative), and that $a(0)\gg1$ (which makes $P_\mathrm{esc}$ nearly
independent of the initial distribution of $x$).

Figure~\ref{fig:snfprob} shows the quasi-static approximation of the
cumulative escape probability $P_\mathrm{esc}$ for $\epsilon=0.1$,
$0.01$ and $0.001$ as grey curves. Superimposed are numerical
observations of the cumulative escape probability $P_\mathrm{esc}$ for
the normal form with drifting parameter $a$
\eqref{eq:snformscal}--\eqref{eq:adriftscal} as black curves. In a
simulation of \eqref{eq:snformscal}--\eqref{eq:adriftscal} one can
approximate the cumulative escape probability $P_\mathrm{esc}$ until
time step $t_n$ with the help of a recursion for the probability
$1-P_\mathrm{esc}$ of \emph{not} escaping
\begin{displaymath}
  1-P_\mathrm{esc}(t_{n+1})=
  \frac{N-N_\mathrm{esc}(t_n)}{N}(1-P_\mathrm{esc}(t_n))
\end{displaymath}
where $N_\mathrm{esc}(t_n)$ is the number of realisations that escape
at time step $t_n$ and $N$ is the overall ensemble size. We
keep the overall ensemble size $N$ (of non-escaped realisations)
constant by re-initialising every escaped realisation to a random
non-escaped instance.

One can see that if the control parameter $a$ drifts slowly
($\epsilon\ll1$) the cumulative escape probability increases sharply
from nearly $0$ to nearly $1$ in a range of $a$ of about length $1$
(for example, between $a=1$ and $a=2$ for $\epsilon=0.001$).

Figure~\ref{fig:percentiles} shows the percentiles of the cumulative
escape probabilities systematically for $\epsilon$ ranging between
$0.001$ and $0.5$. Again, we superimpose black curves showing the
percentiles of the cumulative probability observed during a simulation
of the normal form with drifting control parameter
\eqref{eq:snformscal}--\eqref{eq:adriftscal}.

We draw two conclusions from the results shown in
Figure~\ref{fig:snfprob} and Figure~\ref{fig:percentiles}:
\begin{enumerate}
\item the quasi-static approximation \eqref{eq:probesc} of the
  cumulative escape probability is quantitatively accurate to order
  $\epsilon$. For large $\epsilon$ ($\sim0.1$) the effect of the
  dynamic drifting of the control parameter delays the escape
  slightly. Naturally, this effect is much weaker than the delay in
  exchange of stability observed in slow passages through Hopf or
  Pitchfork bifurcations (see \citep{Kuske1999,Baer1989,Su2004}
  for studies that quantify also the effect of noise).
\item Random disturbances make early escape probable as soon as the
  ratio of drift speed $\epsilon$ and variance of the disturbance
  $\sigma^2$ becomes small in the original parameters of saddle-node
  normal form with drift and noise
  \eqref{eq:gensnform}--\eqref{eq:adrift}. The percentiles in
  Figure~\ref{fig:percentiles} quantify this effect after the axes
  have been scaled using transformation
  \eqref{eq:rescaling}: $\epsilon\mapsto\epsilon\sigma^{-2}$,
  $a\mapsto a\sigma^{-4/3}$.
\end{enumerate}

\paragraph*{Control of accuracy in stochastic simulations}
The ensemble size for the numerical simulation was $N=400$, and the
integration was performed with the Euler-Maruyama scheme (which of
order $1$ for this problem) using stepsize $h=0.01$. Control
calculations using different method parameters give results that are
visually indistinguishable from Figure~\ref{fig:snfprob}. For control
we varied (one-by-one) the stepsize ($h=0.05$), the threshold for a
realisation to count as escaped ($x_\mathrm{th}=-10$), and the
ensemble size ($N=800$). We also varied the re-initialisation
strategy: alternatively we chose the new position of a particle after
escape according to a Gaussian distribution with mean $X_s=\sqrt{a}$
(the bottom of the well) and variance $\sigma^2/\sqrt{a}$. This would
be the stationary distribution obtained for the linearisation of the
re-scaled normal form \eqref{eq:snformscal} in its equilibrium
$X_s$. This alternative re-initialisation works by construction only
for $a>0$.

The Figures~\ref{fig:snfprob} and \ref{fig:percentiles} both show
either long-time limits (such as $k_0(a)$ in Figure~\ref{fig:snfprob})
or cumulative quantities: $K_0(a)$ and $P_\mathrm{esc}$. These
quantities can be approximated more accurately with relatively small
ensemble sizes and simple ensemble integration than, for example,
probability densities. The density of the escape probability is the
time derivative of $P_\mathrm{esc}$, which would be a much ``noisier''
function of time (or $a$) for finite randomly drawn ensembles than
$P_\mathrm{esc}$ unless one applies more sophisticated numerical
methods and uses larger ensemble sizes \citep{Kuske1999,Kuske2000}.

For large $a$ numerical simulations of
the stochastic differential equation \eqref{eq:snformscal} become
inefficient because escape events become extremely rare. However, this
case is covered analytically by classical chemical reaction rate
theory \citep{Hanggi1990}:
\begin{equation}\label{eq:k0alarge}
  k_0(a)\sim\frac{2\sqrt{a}}{\pi}\exp\left(-2\sqrt{a}\right)
  \mbox{\ for $a\gg1$,}
\end{equation}
where the expression on the right-hand side contains the dominant
terms of $a$ in the pre-factor and in the exponent. Thus, replacing
the upper limit of integration $a(0)$ by $+\infty$ in the
approximation \eqref{eq:probesc} of $P_\mathrm{esc}$ gives only a
small change.

\subsection{Extraction of noise-induced escape probability from time series}
\label{sec:extract}
In order to estimate if noise-induced early escape plays a role one
needs to assume that the time series $z_k$ is generated by a system
with a parameter that drifts and approaches a saddle-node
bifurcation. We denote the drifting parameter by $a$ in a manner that
the saddle-node bifurcation is at $a=0$ and equilibria of the
deterministic system with fixed $a$ exist for $a>0$. We introduce the
parameter $Z_0$ as the $z$-value of the equilibrium at the saddle-node
parameter $a=0$. If we assume that the underlying process is subject
to an additive (Gaussian) perturbation of amplitude $\sigma_z$ then
$z_k$ can be written as the measurements of the system
\begin{align}
  \label{eq:snunscaled}
  \d z&=\left[q a(t)-\frac{(z-Z_0)^2}{q}+O((z-Z_0)^3)\right]\d t+
  \sigma_z\d W_t
\end{align}
at times $k\Delta t$. The parameter $q$ in \eqref{eq:snunscaled}
measures the width of the parabola $a(z)$ of equilibria of
\eqref{eq:snunscaled} for zero noise amplitude ($\sigma_z=0$) and
ignoring third-order powers of $z-Z_0$:
\begin{displaymath}
  a(z)=\frac{(z-Z_0)^2}{q^2}\mbox{.}
\end{displaymath}
If $q$ is negative then the stable side of the parabola is below
$Z_0$, if $q$ is positive then the stable side is above $Z_0$. From
the general form \eqref{eq:snunscaled} of a system near a saddle-node
one can obtain the normal form
\begin{equation}\label{eq:ts:snform}
  \d x=\left[a(t)-x^2\right]\d t +\sigma \d W_t
\end{equation}
by ignoring third-order powers of $z-Z_0$ and introducing the
re-scaled variable
\begin{displaymath}
  x=\frac{z-Z_0}{q}\mbox{.}
\end{displaymath}
The noise amplitude scales correspondingly:
\begin{equation}\label{eq:noisescale}
  \sigma=\frac{\sigma_z}{q}\mbox{.}
\end{equation}
We demonstrate how to extract a (crude) estimate of the state and the
parameters of \eqref{eq:snunscaled} from a time series $z_k$ as shown,
for example, for the model output in Figure~\ref{fig:ldr:est}(a,b).
As long as $a(t)$ is large the potential well is deep such that the
Gaussian perturbation $\sigma_z\d W_t$ is unlikely to kick $z(t)$ out
of the well of the stable equilibrium. Choosing $\Delta t$ (the time
spacing of the measurements $z_k$) as our time unit we first revert
the estimates $c_k$ and $\theta_k$ obtained from the degenerate
finger-printing into estimates of the parameters $\kappa(t)$ and
$\sigma(t)$ appearing in the Ornstein-Uhlenbeck process obtained from
linearising \eqref{eq:snunscaled} in its stable equilibrium at
$q\sqrt{\smash[b]{a}}$:
\begin{displaymath}
  \d z=-\kappa(t)z\d t+\sigma_z(t)\d W_t\mbox{.}  
\end{displaymath}
The estimates $\kappa_k$ for $\kappa(k\Delta t)$ and
$\sigma_k$ for $\sigma(k\Delta t)$ are
\begin{align}
  \kappa_k&=-\log c_k\label{eq:kappak}\\
  \intertext{where the estimate $c_k$ is calculated using the fitting
    procedure in Section~\ref{sec:timeseries} and}
  \sigma_{z,k}&=\theta_k\sqrt{\frac{2\kappa_k}{(1-c_k^2)}}\label{eq:sigmak}
\end{align}
(see \citep{Gillespie1992}) where $\theta_k$ is the standard deviation
of the detrended time series from the AR(1) estimate.  The estimate
$\kappa_k$ is an approximation of the linear decay rate
\begin{equation}
  \label{eq:snprop}
  \kappa(k\Delta t)=2\sqrt{a(k\Delta t)}
\end{equation}
as long the time series $z_k$ is near the bottom of the potential
well. Thus, we can estimate $a(k\Delta t)$ for the times $k\Delta t$
as $a_k$:
\begin{displaymath}
  a_k=\frac{\kappa_k^2}{4}\mbox{.}
\end{displaymath}
Consequently, estimates of the parameters $a$ and $\sigma_z$ are
by-products of the AR(1) estimate to obtain $c_k$. The only remaining
unknown quantity that one needs to convert to the normal form
\eqref{eq:ts:snform} is the scaling factor $q$. If we drop third-order
terms of $z-Z_0$, fix $a$ at $a_k$, and consider only the
deterministic part ($\sigma_z=0$) then the equilibrium $Z_k$ of
\eqref{eq:snunscaled} satisfies $\sqrt{\smash[b]{a_k}}q-Z_k+Z_0=0$,
which is equivalent to
\begin{equation}\label{eq:qdet}
  \frac{\kappa_k}{2}q+Z_0=Z_k\mbox{.}
\end{equation}
We see that $q$ is the ratio between the slope of $\kappa_k/2$ (for
which we have an estimate) and the slope of the equilibrium state
which the time series $z_k$ fluctuates around. An estimate for $Z_k$
has also been obtained during the degenerate finger-printing procedure
as the kernel-average of $z_k$ (also named $Z_k$ in
Section~\ref{sec:timeseries}). Thus, an estimate for $q$ is the ratio
between the mean slope of $Z_k$ and the mean slope of $\kappa_k$.

In our analysis of the normal form with drift and noise we considered
the case were the parameter $a$ drifts with a uniform (small) speed
$\epsilon$. This is unrealistic for paleo-climate records and in
models whenever the parameter $a$ is driven by the output of another
subsystem, for example, if $a$ is freshwater forcing as in
\citep{Rahmstorf2000}. In order to predict how the parameter $a$
continues to drift beyond the final sliding window we assume that the
process driving $a$ is stationary. So,
\begin{equation}\label{eq:apredict}
  a(t)=a(0)-\int_0^t\epsilon(s)\d s
\end{equation}
where $\epsilon$ is not constant but a random variable with a
stationary probability distribution. The distribution has been
estimated for the period of time where $\kappa$ is available:
$\epsilon_k$ are the increments between successive $a$ estimates:
\begin{displaymath}
  \epsilon_k=\frac{a_{k}-a_{k+1}}{\Delta t}\mbox{.}
\end{displaymath}
These $\epsilon_k$ form an empirical sample of the distribution of
$\epsilon(s)$ from which we can draw to calculate $a(t)$ using
\eqref{eq:apredict} without estimating any further parameters. If the
mean of $\epsilon$ exists and is bigger than zero then $a(t)$ will reach
$0$ almost surely, resulting in a probability distribution
\begin{displaymath}
  P_a(t)=1-P\left(a(s)>0\phantom{\vert_\vert}\mbox{\ for all $s<t$}\right)\mbox{,}
\end{displaymath}
which is the probability that the random variable $a$ reaches its
critical value $a=0$ before time $t$. The probability
$P_\mathrm{esc}(t)$ of a trajectory escaping before time $t$ can now
be estimated using a direct numerical simulation of
\begin{displaymath}
  \begin{split}
   \d x&=\left[a-x^2\right]\d t +\sigma \d W_t\\
    \d a&=-\epsilon(t)\d t\mbox{.}
  \end{split}
\end{displaymath}
Alternatively, if the mean of $\epsilon$ is sufficiently small one can
use the quasi-static approximation
\begin{displaymath}
  P_\mathrm{esc}(t)=1-E\left[\exp\left(\int_0^tk_0(a(s))\d s\right)\right]
\end{displaymath}
where the escape rate $k_0(a)$ (which exists only for positive $a$) is
given in Figure~\ref{fig:snfprob} and $E[x]$ is the mean of the random
variable $x$.

\begin{figure}[t!]
  \centering
  \includegraphics[width=\textwidth]{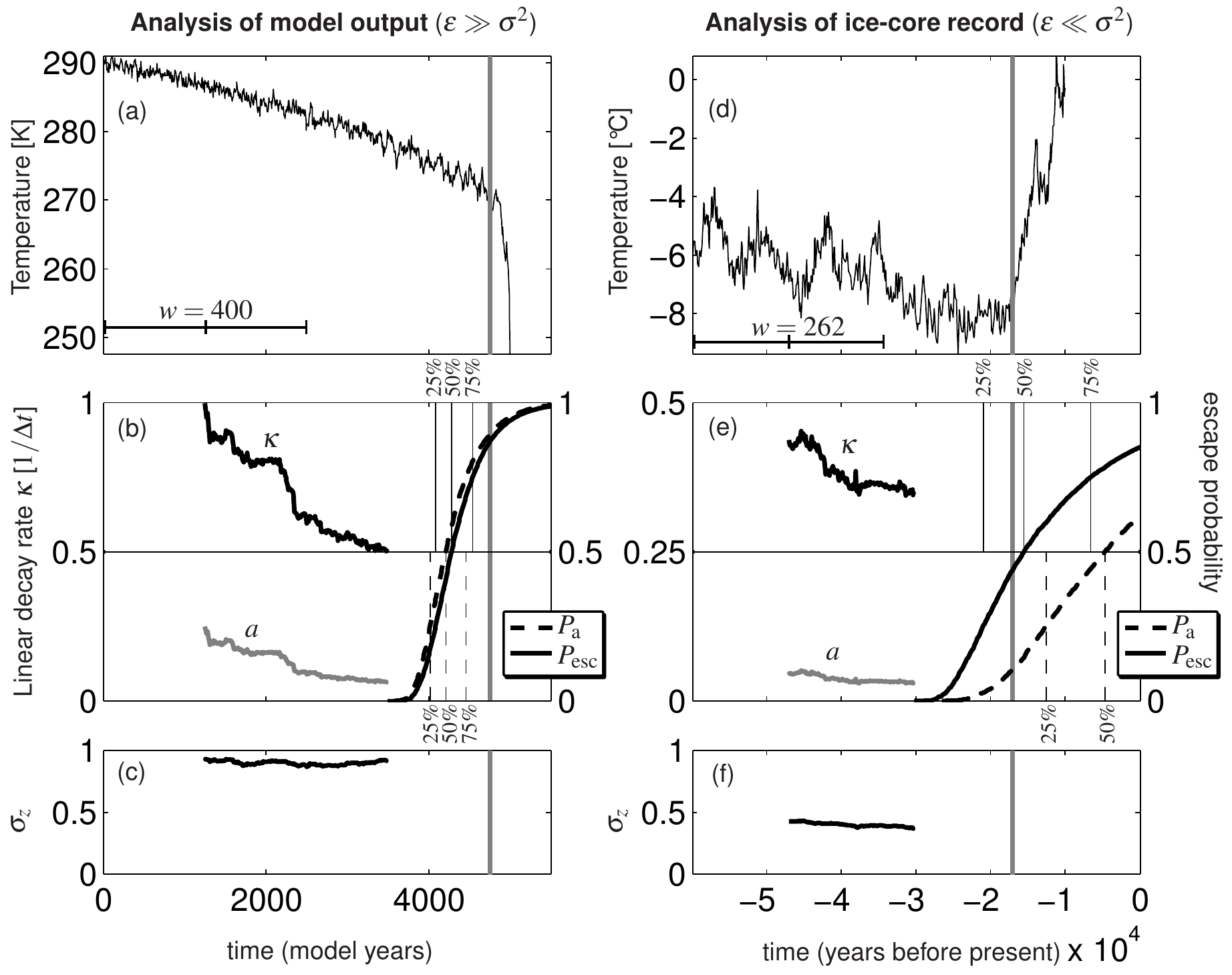}
  \caption{Prediction of saddle-node bifurcation and early escape for
    two time series shown in \citep{Dakos2008}. Diagram (a) shows the
    stochastic model output also shown in
    Figure~\ref{fig:ldr:est}(a). All quantities have been obtained
    only from data prior to the cut-off time indicated by the grey
    vertical line in diagram (a). Diagram (b) shows the linear decay
    rate $\kappa_k$ (as extracted by degenerate finger-printing with
    window size $w$, compare to Figure~\ref{fig:ldr:est}(b)) and the
    corresponding normal form parameter $a$. It also contains the
    cumulative probability functions for reaching the critical value
    $a=0$ and for escape, $P_a$ and $P_\mathrm{esc}$ on the right
    axis. The graph in diagram (c) is the estimated noise amplitude
    $\sigma$ extracted from the time series. The diagrams (d-f)
    discuss the time series shown in Figure~\ref{fig:ldr:est}(c),
    which is a snapshot from the paleo-climate record of
    \citet{Petit1999} at the end of the last glaciation using the
    identical procedure to the one applied to the time series in
    diagram (a).}
  \label{fig:snf:timeseries}
\end{figure}
\begin{table}[ht!]
  \centering
  \begin{tabular}[t]{ccc}\hline\noalign{\smallskip}
    Parameter & \qquad Model output Figure~\ref{fig:snf:timeseries}(a--c)\qquad &
    \qquad Ice-core record Figure~\ref{fig:snf:timeseries}(d--f)
    \\\noalign{\smallskip}\hline\noalign{\smallskip}
    mean $\epsilon$ & $2.2\times10^{-4}$ & $5.2\times10^{-6}$\\
    stdev $\epsilon$ & $1.1\times10^{-3}$ & $5.9\times10^{-5}$\\
    $q$ & $621$ & $13.5$\\
    $\sigma$ & $1.4\times10^{-3}$ & $3.0\times10^{-2}$\\
    mean $\epsilon\sigma^{-2}$ & $105$ & $5.9\times10^{-3}$\\\hline
  \end{tabular}
  \caption{Numerical values of estimates of 
    parameters that are not visible in 
    Figure~\ref{fig:snf:timeseries}.}
  \label{tab:numbers}
\end{table}
Figure~\ref{fig:snf:timeseries}(b) and (c) show the quantities
$\kappa_k$, $a_k$ and $\sigma_{z,k}$ for the time series $z_k$ of the
stochastic model output in Figure~\ref{fig:snf:timeseries}(a) (the
same series $z_k$ as in
Figure~\ref{fig:ldr:est}(a)). Figure~\ref{fig:snf:timeseries}(b) also
shows the cumulative probability function $P_a(t)$ for $a$ crossing
the critical value $0$ before $t$ (dashed curve) and its quartiles
(dashed threshold lines), and the probability $P_\mathrm{esc}(t)$ for
escape before time $t$ (solid curve) and its quartiles (solid
threshold lines). The threshold for a trajectory counting as escaped
was set at $x_\mathrm{th}=-5$.  Table~\ref{tab:numbers} lists the
values of $q$ as determined by relation \eqref{eq:qdet} and
$\sigma=\sigma_z/q$. Also, the table lists the empirical mean value of
the random variable $\epsilon$ and its empirical standard deviation
(the distribution of $\epsilon$ is only moderately skewed). The mean
of $\epsilon/\sigma^2$ indicates if the parameter drift is rapid or
slow in rescaled normal form
coordinates. Figure~\ref{fig:snf:timeseries}(e) and (f) show the
quantities $\kappa$, $a$, $P_\mathrm{esc}$, $P_a$ and $\sigma_z$ for
the ice-core record in Figure~\ref{fig:snf:timeseries}(d) (the same
series $z_k$ as in Figure~\ref{fig:ldr:est}(c))

We observe that the parameter drift in the model output (see
Figure~\ref{fig:snf:timeseries}(b)) is fast compared to the noise
level such that early escape plays no role. Table~\ref{tab:numbers}
also shows that the ratio between drift and variance is large, which
confirms the visual impression that drift dominates the noise. The
probability distribution for escape is even shifted to the right such
that trajectories are expected to escape \emph{later} than the
drifting system parameter $a$ reaches its critical value.

We note that the actual escape of the observed instance from
Figure~\ref{fig:snf:timeseries}(a) occurs relatively late (at the
$90\%$ percentile). This shifts down (to below $75\%$) for shorter
sliding windows in the degenerate finger-printing procedure. Both
probability distributions are relatively symmetric and concentrated in
a range of approximately $1\,000$ model years.

In contrast, the time series of the ice-core record, shown in
Figure~\ref{fig:snf:timeseries}(d), has a slowly drifting parameter
compared to the noise level (note that one cannot be certain that the
underlying mechanism for the apparent tipping is indeed a passage
through a saddle-node). Consequently, trajectories of the estimated
saddle-node normal form system escape significantly earlier than $a$
reaches its critical value: the cumulative probability
$P_\mathrm{esc}$ is shifted to the left of $P_a$. For example, the
median ($50\%$) time for escape is $1.1\times10^{4}$ years before the
median time for reaching the critical value $a=0$. We notice that both
distributions, $P_a$ and $P_\mathrm{esc}$ are highly skewed, having a
long tail at the right end. This makes the expected escape time (mean
escape time) as a point estimate sensitive to small perturbations such
as measurement uncertainty or a different choice of method parameters
in the degenerate finger-printing procedure whereas the median times
are comparatively robust. A systematic analysis, including the
dependence on the cut-off time (the grey vertical line in all panels
of Figure~\ref{fig:snf:timeseries}), is a topic for future work.



\section{Conclusion}
\label{sec:conc}

Methods to identify incipient climate tipping using time series
analysis have recently been developed by \citet{Held2004} and
\citet{Livina2007}. These methods have been tested on model outputs
and paleo-climate data by \citet{Livina2007} and
\citet{Dakos2008}. Time series analysis should be seen as a complement
to the huge modelling efforts that are invested into the analysis and
prediction of climate changes. While \citet{Dakos2008} could
demonstrate that the characteristic quantity extracted from the time
series, the propagator, indeed increases (as it should according to
bifurcation theory) using this quantity for prediction is much more
challenging. The extraction of the propagator from a time series makes
assumptions that are, for archaeological records, difficult to check:
for example, separation of time scales between parameter drift, decay
of critical mode and decay of stable modes, or the nearly
constant-speed approach of the underlying control parameter toward its
critical value. 

We studied another source of uncertainty, which exists even if the
underlying deterministic dynamics is drifting with constant speed
close to a saddle-node bifurcation and the estimated linear decay rate
is accurate. Namely, the escape of the dynamics from the potential
well around the stable equilibrium can be either premature or delayed,
depending on the ratio between parameter drift speed and the amplitude
of the random disturbances. We derived an approximate (semi-analytic)
formula that is valid in the quasi-static limit (that is, the
parameter drifts sufficiently slowly compared to the noise
amplitude). We found that the early escape effect vanishes if the
drift is more rapid (drift speed $\epsilon\approx0.3$ in normal form
scaling). 

We also demonstrated how one can estimate this effect from time series
data, using what we might call a `fold for incipient tipping' method
(FIT). This was tested for two examples, the output of a stochastic
model (a case of rapid drift) and an ice-core record (a case of slow
drift, or large noise amplitude). We plan to study the consistency of
the proposed escape prediction for time series more systematically in
the future. This should be done by generating time series instances
from saddle-node normal forms with drifting parameter and noise,
predicting escape distributions from these time series, and then
comparing the predictions to the original distributions shown in
Figure~\ref{fig:percentiles}.

Our estimates for early escape from a potential well are all stated in
probabilistic terms, which is appropriate if one treats the
disturbances as random. Normal form based estimates such as ours may
also be useful when studying tipping in climate models because running
simulations for large ensembles of realizations in sophisticated
climate models is expensive. If the modelled scenario is close to
tipping and deterministic calculations reveal the presence of a
saddle-node bifurcation then the diagrams in Figure~\ref{fig:snfprob} and
Figure~\ref{fig:percentiles} help to give estimates for the escape
probability depending on the speed of parameter drift and the
amplitude of the disturbances.

Finally, we note that our scenario corresponds to the escape from a
potential well of an overdamped particle (this appears to be the
relevant case, for example, in ocean circulation models, see
\citep{Dijkstra2004}). The weakly damped case,
\begin{displaymath}
  \frac{\d^2x}{\d t^2}+\gamma \frac{\d x}{\d t}=-U(x)+\sigma\xi(t)\mbox{,}
\end{displaymath}
where $\gamma$ is a small coefficient determining the amount of
damping ($0<\gamma\ll1$), $U$ is the potential and $\xi$ is the
disturbance (for example, a random Gaussian noise), has also been
covered by reaction rate theory (the classical theory treats the case
where the noise amplitude $\sigma$ is much smaller than the potential
barrier, see \citep{Hanggi1990} for a review).

%

%

\end{document}